\documentclass[11pt]{article}
\usepackage{amsmath,amssymb,amsfonts,amsthm}

\newcommand{\Ann}{\mbox{Ann}\,}

\newcommand{\Ass}{\mbox{Ass}\,}

\newcommand{\gr}{\mbox{girth}\,}

\renewcommand{\d}{\mbox{d}\,}

\newcommand{\rad}{\mbox{rad}}

\newcommand{\E}{\mbox{E}}

\newcommand{\diam}{\mbox{diam}\,}

\newcommand{\clique}{\mbox{clique}\,}

\newcommand{\G}{\mbox{G}}
\newcommand{\V}{\mbox{V}}
\newcommand{\Z}{\mbox{Z}}
\newcommand{\e}{\mbox{e}}
\newcommand{\C}{\mbox{C}}

\newcommand{\fp}{\frak{p}}
\newcommand{\fq}{\frak{q}}

\begin{document}

\date{}

\title{\bf Complete $r$-partite zero-divisor graphs and coloring of commutative semigroups
\footnotetext{* Corresponding author. Department of Math., Univ. of
Tehran, P.O. Box 13145-448, Tehran, Iran}
}

\author{H. R. Maimani $^{\it ab}$, M. Mogharrab $^{\it c}$ and S. Yassemi $^{\it
ab*}$\\
{\small\it $(a)$ Department of Mathematics, University of Tehran,
Tehran, Iran}\\
{\small\it $(b)$ Institute for Studies in Theoretical Physics and
Mathematics, Tehran, Iran}\\
{\small\it $(c)$ Department of Mathematics, Persian Golf University,
Bushehr, Iran}}


\maketitle

\begin{abstract}

\noindent For a commutative semigroup $S$ with $0$, the
zero-divisor graph of $S$ denoted by $\Gamma(S)$ is the graph
whose vertices are nonzero zero-divisor of $S$, and two vertices
$x$, $y$ are adjacent in case $xy=0$ in $S$. In this paper we
study the case where the graph $\Gamma(S)$ is complete $r$-partite
for a positive integer $r$. Also we study the commutative
semigroups which are finitely colorable.

\end{abstract}


\vspace{.3in}

\noindent{\bf 1. Introduction}

\vspace{.2in}

In [{\bf B}] Beck introduced the concept of a zero-divisor graph
$\G(R)$ of a commutative ring $R$. However, he lets all elements
of $R$ be vertices of the graph and his work was mostly concerned
with coloring of rings. Later, Anderson and Livingston in [{\bf
AL}] studied the subgraph $\Gamma(R)$ of $\G(R)$ whose vertices
are the nonzero zero-divisors of $R$. The zero-divisor graph of a
commutative ring has been studied extensively by several authors,
e.g. [{\bf AFL}], [{\bf ALS}], [{\bf LS}], and [{\bf AMY}].

For the sake of completeness, we state some definitions and notions
used throughout to keep this paper as self contained as possible.

\indent For a graph $G$, the {\it degree} of a vertex $v$ in $G$
is the number of edges of $G$ incident with $v$. For a nontrivial
connected graph $G$ and a pair $u,v$ of vertices of $G$, the
distance $d(u,v)$ between $u$ and $v$ is the length of shortest
path from $u$ to $v$ in $G$. The {\it eccentricity} $\e(v)$ of a
vertex $v$ in graph $G$ is the distance from $v$ to a vertex
farthest from $v$, that is,
$$\e(v)=\max\{\d(x,v)|x\in\V(G)\}.$$
The {\it radius} $\rad(G)$ of a connected graph is defined as
$$\rad(G)=\min\{\e(v)|v\in\V(G)\},$$
and the {\it diameter} $\diam(G)$ of a connected graph $G$ is
defined as
$$\diam(G)=\max\{\e(v)|v\in\V(G)\}.$$
It is known that $$\rad(G)\le\diam(G)\le 2\,\rad(G).$$ An {\it
$r$-partite} graph is one whose vertex set can be partitioned into
$r$ subsets so that no edge has both ends in any one subset. A {\it
complete $r$-partite} graph is one in which each vertex is joined to
every vertex that is not in the same subset. The {\it complete
bipartite} (i.e., $2$-partite) graph is denoted by $K_{m,n}$ where
the set of partition has sizes $m$ and $n$. A graph in which each
pair of distinct vertices is joined by an edge is called a {\it
complete} graph. We use $K_n$ for the complete graph with $n$
vertices. The {\it girth} of a graph $G$ is the length of a shortest
cycle in $G$ and is denoted by $\gr(G)$. We define a {\it coloring}
of a graph $G$ to be an assignment of colors (elements of some set)
to the vertices of $G$, one color to each vertex, so that adjacent
vertices are assigned distinct colors. If $n$ colors are used, then
the coloring is referred to as an {\it $n$-coloring}. If there
exists an $n$-coloring of a graph $G$, then $G$ is called
$n$-colorable. The minimum $n$ for which a graph $G$ is
$n$-colorable is called the {\it chromatic number} of $G$, and is
denoted by $\chi(G)$. A {\it clique} of a graph is a maximal
complete subgraph and the number of vertices in the largest clique
of graph G, denoted by $\omega(G)$, is called the {\it clique
number} of $G$. Obviously $\chi(G)\ge\omega(G)$ for general graph
$G$ (see [{\bf CO}, page 289]).

Suppose that $S$ ba a semigroup. A non-empty subset $I$ of $S$ is
called ideal if $xS\subseteq I$ for any$x\in I$. An ideal $\fp$ of a
semigroup is called a prime ideal of $S$ if $xSy\subseteq \fp$
implies $x\in \fp$ or $y\in \fp$. Now let $S$ be a commutative
semigroup with $0$ and let $\Z(S)$ be its set of zero-divisors.
According to [{\bf DMS}], the zero-divisor graph, $\Gamma(S)$, is an
undirected graph with vertices ${\Z(S)^*}=\Z(S)\setminus \{0\}$, the
set of nonzero zero-divisors of $S$, where for distinct $x,y \in
{\Z(S)^*}$, the vertices $x$ and $y$ are adjacent if and only if
$xy=0$. In [{\bf DMS}] DeMeyer, McKenzie, and Schneider observe that
$\Gamma(S)$ (as in the ring case) is always connected, and the
diameter of $\Gamma(S)\le 3$. If $\Gamma(S)$ has a cycle then
$\gr(\Gamma(S))\le 4$. They also show that the number of minimal
ideals of $S$ gives a lower bound to the clique number of $S$. In
[{\bf ZW}] Zue and Wu studied a graph $\bar{\Gamma}(S)$ where the
vertex set of this graph is $\Z(S)^*$ and for distinct elements
$x,y\in\Z(S)^*$, if $xSy=0$, then there is an edge connecting $x$
and $y$. Note that $\Gamma(S)$ is a subgraph of $\bar{\Gamma}(S)$.
Recently, F. DeMeyer and L. DeMeyer studied further the graph
$\Gamma(S)$ and its extension to a simplicial complex, cf. [{\bf
DD}]. Clearly for any prime ideal $\fp$ if $x$ and $y$ are adjacent
in $\Gamma(S)$, then $x\in\fp$ or $y\in\fp$. So for every prime
ideal $\fp$ and every edge $e$, one of the end points of $e$ belongs
to $\fp$

 One may
address three major problems in this area: characterization of the
resulting graphs, characterization of the semigroups with
isomorphic graphs and realization of the connections between the
structures of a semigroup and the corresponding graph. In this
paper we focus on the third problem.

The organization of this paper is as follows:

In Section 2, among the other things, it is shown that for a
reduced commutative semigroup $S$, if $\Z(S)\neq\{0\}$ and every
vertex of $\Gamma(S)$ has finite degree, then $\Z(S)$ is finite,
see Proposition 2.9. It is also shown that if the set of
associated primes of $S$, $\Ass(S)$, has more than two elements
then the girth of $\Gamma(S)$ (i.e. the length of the shortest
cycle in $\Gamma(S)$) is three.

In Section 3, we study the semigroups whose zero divisor graphs are
complete $r$-partite. It is shown that for a reduced commutative
semigroup $S$ if $\Gamma(S)$ is a complete $r$-partite graph, with
parts $V_1, V_2,..., V_r$, then $V_t\cup\{0\}$ is an ideal and
$\fp_t=Z(S)\setminus V_t$ is a prime ideal for any $1\leq t\leq r$.

In Section 4, we study the semigroups of finite chromatic number.
 We show that for
a commutative semigroup $S$ the following conditions are
equivalent: (1) $\chi(S)<\infty$, (2) $\omega(S)<\infty$, and (3)
the zero ideal is a finite intersection of prime ideals, where
$\chi(S)=\chi(\Gamma(S))$ and $\omega(S)=\omega(\Gamma(S))$, see
Theorem 4.1. As a corollary we show that $\chi(S)=\omega(S)=n$ if
$S$ is a reduced semigroup and $0=\cap_{i=1}^n\fp_i$ is a minimal
prime decomposition of $0$ (i.e. for any $i\neq j$,
$\fp_i\neq\fp_j$ and for any $1\le t\le n$, $0\neq\cap_{i\neq
t}\fp_i$). In addition, it is shown that for $n\le 2$, $\chi(S)=n$
if and only if $\omega(S)=n$. It is shown that this result is not
valid for $n=3$. We give a finite commutative semigroup $S$ with
$\chi(S)=4$ and $\omega(S)=3$.

We follow standard notation and terminology from graph theory [{\bf
CO}] and semigroup theory [{\bf H}].

\vspace{.3in}

\noindent{\bf 2. Some special ideals and girth of $\Gamma(S)$}

\vspace{.2in}

Let $S$ be a commutative semigroup with $0$. It is known that the
following hold:

\begin{itemize}

\item[(a)] $Z(S)$ is an ideal of $S$;

\item[(b)] $S'=S \setminus Z(S)$ and $S'\cup {0}$ are
subsemigroup of $S$ with no nonzero zero-divisors.

\end{itemize}

Let $T$ be a non-empty set of vertices of the graph $G$. The
subgraph induced by $T$ is the greatest subgraph of $G$ with
vertex set $T$, and is denoted by $G[T]$, that is, $G[T]$ contains
precisely those edges of $G$ joining two vertices of $T$.

The following result gives a graph property of the subgraph of
$\Gamma(S)$ which consists of the nonzero nilpotent elements of
$S$.

\vspace{.1in}

\noindent {\bf Proposition 2.1.} Let $N$ be the set of nilpotent
elements of $S$. If $N^*=N\setminus\{0\}$ is a non-empty set, then
$\Gamma(S)[N^*]$ is a connected subgraph of $\Gamma(S)$ of
diameter at most $2$.

\vspace{.1in}

\indent {\it Proof.} Since $N$ is a semigroup we have that
$\Gamma(N)=\Gamma(S)[N^*]$ is connected, see [{\bf DMS}, Theorem
1.2]. In addition, $N$ is nilpotent semigroup and so $\diam
\Gamma(N)\le 2$, see [{\bf DD},Theorem 5].\hfill$\square$

\vspace{.2in}

 The {\it distance $d(v)$} of
a vertex $v$ in a connected graph $G$ is the sum of the distances
$v$ to each vertex of $G$. The {\it median $M(G)$} of a graph $G$ is
the subgraph induced by the set of vertices having minimum distance.

 Let $G$ be a connected graph, and $T\subseteq \V(G)$. We say $T$ is
 a {\it cut vertex set} if $G\setminus T$ is disconnected. Also the cut vertex set $T$ is
 called a minimal cut vertex set for $G$ if no proper subset of
 $T$ is a cut vertex set. In addition, if $T=\{x\}$, then $x$ is
 called a {\it cut vertex}.

\vspace{.1in}

\noindent {\bf Theorem 2.2.} The set of vertices of $M(\Gamma(S))
\bigcup \{0\}$ is an ideal of S. In addition, if $T$ is a minimal
cut vertex set of $\Gamma(S)$, then $T\cup\{0\}$ is an ideal of
$S$.

\vspace{.1in}

\noindent {\it Proof.} Let $x$ be a vertex of $M(\Gamma(S))$ and
$y\in S$. Suppose that $xy\neq 0$. Let $z$ be a vertex of
$\Gamma(S)$ and $\d(x,z)=t$. Then there is a shortest path from $x$
to $z$ of length $t$,
\begin{center}
$x$---$x_1$---$x_2$---$\cdots$---$x_{t-1}$---$z$
\end{center}
and so
\begin{center}
$xy$---$x_1$---$x_2$---$\cdots$---$x_{t-1}$---$z$,
\end{center}
is a walk of length $t$ from $xy$ to $z$. Thus $\d(xy,z)\le\d(x,z)$.
Since $\d(r,r)=0$, we have the following (in)equalities:
$$\d(xy)=\sum_{z\in\V(\Gamma(S))}\d(xy,z)\le\sum_{z\in\V(\Gamma(S))}\d(x,z)=\d(x).$$
Since $x\in M(\Gamma(S))$, we have $\d(xy)=\d(x)$, and hence $xy$
belongs to the vertex set of $M(\Gamma(S))$.

Now let $T$ be a minimal cut vertex set of $\Gamma(S)$, and $x\in
T$, $r\in S$. Since $T\setminus\{x\}$ is not a cut vertex of
$\Gamma(S)$, there exist two vertices $z,y$ of the graph $\Gamma(S)$
such that $y$---$x$---$z$ is a path in $\Gamma(S)$, and $y,z$ belong
to two distinct connected components of $\Gamma(S)\setminus T$. Now
if $rx\neq 0$, and $rx\notin T$, then $rx$ is a vertex of
$\Gamma(S)\setminus T$. Therefore we have the following path in
$\Gamma(S)\setminus T$;
\begin{center}
$y$---$rx$---$z$,
\end{center}
which is a contradiction. Thus $rx\in T\cup\{0\}$ and so
$T\cup\{0\}$ is an ideal of $S$.\hfill$\square$

\vspace{.2in}

The techniques of the proof of Theorem 2.2 can be applied to obtain
the following result.

\vspace{.1in}

\noindent{\bf Corollary 2.3.} Let $x$ be a cut vertex of
$\Gamma(S)$. Then $\{0,x\}$ is an ideal of $S$. In this case either
$x$ is adjacent to every vertex of $\Gamma(S)$ or $x\in Sx$.

\vspace{.2in}

The {\it center} $\C(G)$ of a connected graph $G$ is the subgraph
induced by the vertices of $G$ with eccentricity equal the radius of
$G$.

\vspace{.1in}

\noindent{\bf Theorem 2.4.} For the semigroup $S$, the set
$\V(\C(\Gamma(S)))\cup\{0\}$ is an ideal of $S$.

\vspace{.1in}

\noindent{\it Proof.} Let $x\in\V(\C(\Gamma(S)))$, and $r\in S$.
Suppose that $rx\neq 0$. Then
$$\e(rx)=\max\{\d(u,rx)|u\in\V(G)\}\le\max\{\d(u,x)|u\in\V(G)\}=\e(x).$$
Thus $\e(rx)=\e(x)$, and so
$rx\in\V(\C(\Gamma(S)))\cup\{0\}$.\hfill$\square$

\vspace{.2in}

A subgraph $H$ of a graph $G$ is a {\it spanning} subgraph of $G$
if $\V(H)=\V(G)$. If $U$ is a set of edges of a graph $G$, then
$G\setminus U$ is the spanning subgraph of $G$ obtained by
deleting the edges in $U$ from $\E(G)$. A subset $U$ of the edge
set of a connected graph $G$ is an {\it edge cutset} of $G$ if
$G\setminus U$ is disconnected. An edge cutset of $G$ is {\it
minimal} if no proper subset of $U$ is edge cutset. If $e$ is an
edge of $G$, such that $G\setminus\{e\}$ is disconnected, then $e$
is called a {\it bridge}. Note that if $U$ is a minimal edge
cutset, then $G\setminus U$ has exactly two connected components.

\vspace{.1in}

\noindent{\bf Theorem 2.5.} Let $e=xy$ be a bridge of $\Gamma(S)$
such that the two connected components $G_1$, $G_2$ of
$\Gamma(S)\setminus\{e\}$ have at least two vertices. Then
$Sx=\{0,x\}$ and $Sy=\{0,y\}$ are two minimal ideals of $S$. Also
if $G_1$ or $G_2$ has only one vertex (i.e. $\deg x=1$ or $\deg
y=1$), then $\{0,x,y\}$ is an ideal.

\vspace{.1in}

\noindent{\it Proof.} Since $G_1$ and $G_2$ have at least two
vertices, there exists vertices $g_1$ and $g_2$ of $\Gamma(S)$ with
$g_1\in\V(G_1)$, $g_2\in\V(G_2)$, and $x$ adjacent to $g_1$ (in
$G_1$) and $y$ adjacent to $g_2$ (in $G_2$). Suppose that $r\in S$
and $rx\neq 0$. Then $rx\in\Z(S)$. If $rx\in G_2$, then $rx$ is
adjacent to $g_1$ in $\Gamma(S)\setminus\{e\}$, which is a
contradiction. Therefore $rx\in G_1$. We claim that $rx=x$. In the
other case $rx$ is adjacent to $y$ in $\Gamma(S)\setminus\{e\}$,
which is a contradiction. Since $g_2x\neq 0$ we have that $g_2x=x$
and so $Sx=\{0,x\}$ is a minimal ideal of $S$. Similarly
$Sy=\{0,y\}$ is a minimal ideal of $S$. The last part follows by a
similar argument.\hfill$\square$

\vspace{.2in}

The techniques of the proof of Theorem 2.5 can be applied to obtain
the following result.

\vspace{.1in}

\noindent{\bf Corollary 2.6.} Let $T$ be the minimal edge cutset of
$\Gamma(S)$, and $G_1$, $G_2$ are two parts of $G\setminus T$. Then
the following hold.

\begin{itemize}

\item[(a)] For any $i=1,2$, $(\V(G_i)\cap\V(T))\cup\{0\}$ is ideal of $S$ provided
$G_i$ has at least two vertices.

\item[(b)] $\V(T)\cup\{0\}$ is an ideal if $G_1$ or $G_2$ has only one
vertex.

\end{itemize}

\vspace{.2in}

A semigroup is called reduced if for any $x\in S$, $x^n=0$ implies
$x=0$. We define the annihilator as a non-zero ideal of the form
$\Ann(x)$ for some $x\in S$.

\vspace{.1in}

\noindent{\bf Proposition 2.7.} Let $S$ be a reduced semigroup
which $\Gamma(S)$ does not contain an infinite clique. Then $S$
satisfies the a.c.c on annihilators.

\vspace{.1in}

\noindent{\it Proof.} Suppose that $\Ann x_1<\Ann x_2<\cdots$ be
an increasing chain of ideals. For each $i\ge 2$, choose
$a_i\in\Ann x_i\setminus\Ann x_{i-1}$. Then each $y_n=x_{n-1}a_n$
is nonzero, for $n=2,3,\cdots$. Also $y_iy_j=0$ for any $i\neq j$.
Since $S$ is a reduced semigroup, we have $y_i\neq y_j$ when
$i\neq j$. Therefore we have an infinite clique in $S$. This is a
contradiction and so the assertion holds.\hfill$\square$

\vspace{.2in}

\noindent{\bf Lemma 2.8.} Let $S$ be a commutative semigroup and let
$\Ann a$ be a maximal element of $\{\Ann x:0\neq x\in S\}$. Then
$\Ann a$ is a prime ideal.

\vspace{.1in}

\noindent{\it Proof.} Let $xSy\subseteq\Ann a$, and $x, y\notin\Ann
a$. Then $xxy\in\Ann a$, and so $x^2ya=0$. Since $ya\neq 0$ and
$\Ann a\subset\Ann ya$, we have $\Ann a=\Ann ya$. Thus $x^2\in\Ann
a$ and hence $x\in\Ann xa=\Ann a$. This is a
contradiction.\hfill$\square$

\vspace{.2in}

Recall that the set of associated primes of a commutative
semigroup $S$ is denoted by $\Ass(S)$ and it is the set of prime
ideals $\fp$ of $S$ such that there exists $x\in S$ with
$\fp=\Ann(x)$.  The next result gives some information of
$\Gamma(S)$.

\vspace{.1in}

\noindent{\bf Proposition 2.9.} Let $S$ be a commutative semigroup.
Then the following hold:

\begin{itemize}

\item[(a)] If $|\Ass(S)|\ge 2$ and $\fp=\Ann(x)$, $\fq=\Ann(y)$ are two distinct elements of $\Ass(S)$, then $xy=0$.

\item[(b)] If $|\Ass(S)|\ge 3$, then $\gr(\Gamma(S))=3$.

\item[(c)] If $|\Ass(S)|\ge 5$, then $\Gamma(S)$ is not planar
(A graph $G$ is planar if it can be drawn in the plane in such a way
that no two edges meet except at vertex with which they are both
incident).

\end{itemize}

\vspace{.1in}

\noindent{\it Proof.} (a). We can assume that there exists
$r\in\fp\setminus\fq$. Then $rx=0$ and so $rSx=0\in\fq$. Since $\fq$
is a prime ideal, $x\in\fq$ and hence $xy=0$.

(b). Let $\fp_1=\Ann(x_1)$, $\fp_2=\Ann(x_2)$, and $\fp_3=\Ann(x_3)$
belong to $\Ass(S)$. Then $x_1$---$x_2$---$x_3$---$x_1$ is a cycle
of length 3.

(c). Since $|\Ass(S)|\ge 5$, $K_5$ is a subgraph of $\Gamma(S)$, and
hence by Kuratowski's Theorem $\Gamma(S)$ is not planar .
\hfill$\square$

\vspace{.3in}

\noindent{\bf 3. Complete $r$-partite graph}

\vspace{.2in}

Let $R$ be an infinite ring and let the zero-divisor graph of $R$,
$\Gamma(R)$, be a complete $r$-partite with parts $V_1, V_2,
\cdots ,V_r$ and $r\ge 3$. In [{\bf AMY}, Theorem 3.5] it is shown
that for any integer $1\le t\le r$ and for any $x\in V_t$,
$Rx\subseteq V_t\cup\{0\}$, and $\cup_{i\neq t}V_i\cup\{0\}$ is a
prime ideal. In the following we give a semigroup version of this
result.

\vspace{.1in}

\noindent{\bf Theorem 3.1.} Let $S$ be a reduced commutative
semigroup and let $\Gamma(S)$ be a complete $r$-partite graph with
parts $V_1, V_2,..., V_r$. Then $V_t\cup\{0\}$ is an ideal and
$\fp_t=Z(S)\setminus V_t$ is a prime ideal for any $1\leq t\leq r$.

\vspace{.1in}

\noindent {\it Proof.} For an arbitrary integer $1\le t\le r$
choose $x\in V_t$ and $r\in S$ such that $rx\neq 0$. For any
$i\neq t$, there exists $x_i\in V_i$ with $x_ix=0$. Then
$x_i(rx)=0$. Since $S$ is reduced we have $x_i\neq rx$ for all
$i\neq t$ and hence $rx\in V_t$. Therefore $V_t\cup\{0\}$ is an
ideal. By the same argument $\fp_t$ is an ideal. Now suppose that
$xSy\subseteq \fp_t$, and $s_1\in V_t$. Then $xs_1y\in\fp_t$, and
so $xs_1y=0$. If $xs_1\neq 0$, then $xs_1$---$y$ and $y\notin
V_t$; otherwise $x$---$s_1$ and $x\notin V_t$. Therefore
$x\in\fp_t$ or $s_1y=0$. That implies $x\in\fp_t$ or $y\in\fp_t$.
Thus $\fp_t$ is a prime ideal.\hfill$\square$

\vspace{.2in}

\noindent{\bf Remark 3.2.} (a) It is easy to see that we can
replace the condition ``reduced'' with ``for every $x\in
S\setminus {0}$, $x^2\neq 0$'' in the Theorem 3.1.

(b) in Theorem 3.1 if $\Gamma(S)$ is bipartite (i.e. $r=2$), then
$\Gamma(S)$ is guaranteed to be a complete bipartite graph.

\vspace{.2in}

\noindent{\bf Corollary 3.3.} If $\Ass(S)=\{\fp_1 ,\fp_2\}$,
$|\fp_i|\ge 3$ for $i=1,2$ and $\fp_1\cap\fp_2=\{0\}$, then
$\gr(\Gamma(S))=4$.

\vspace{.1in}

\noindent{\it Proof.} By Theorem 3.1 and Remark 3.2(b), since
$\Gamma(S)$ is complete bipartite, we have that
$\gr(\Gamma(S))=4$.\hfill$\square$

\vspace{.2in}

The following examples show that the condition ``reduced'' is not
redundant in the Theorem 3.1.

\vspace{.1in}

\noindent{\bf Example 3.4.} Let $S=\{0,a,b,c,d\}$ with
$b^2=ab=bc=cd=0$, $ac=c^2=a^2=c$, $d^2=d$, $ad=bd=b$. Then
$\Gamma(S)$ is a bipartite graph as shown in the following
diagram:
\begin{center}
$a$---$b$---$c$---$d$,
\end{center}
where the two parts of $\Gamma(S)$ are $V_1=\{a,c\}$ and
$V_2=\{b,d\}$. It is easy to see that $\{a,c,0\}$ is not an ideal.

\vspace{.2in}

\noindent{\bf Example 3.5.} Let $S=\{0,x,y,z\}$ with $z^2=yz=xz=0$,
$yx=x$, $x^2=x$, $y^2=y$. Then $\Gamma(S)$ is a bipartite graph. In
this case $\{0,x,y\}$ is an ideal but it is not a prime ideal.

\vspace{.2in}

The condition ``reduced'' is not redundant in the statement of
Theorem 3.1. However, it may be replaced by the condition
``$|V_i|>1$ for all $i$'' as we outline below.

\vspace{.1in}

\noindent{\bf Theorem 3.6.} Suppose that $\Gamma(S)$ is complete
$r$-partite graph with parts $V_1, V_2, \cdots ,V_r$ such that for
any $i$, $|V_i|>1$. Then $S$ is reduced.

\vspace{.1in}

\noindent {\it Proof.} Let $x\in V_i$ and $r\in S$ such that
$rx\neq 0$. Since for any $i\neq t$ $|V_i|>1$, there exists
$x_i\in V_i$ such that $rx\neq x_i$ but $rxx_i=0$. Thus $rx\in
V_t$. By the same argument as Theorem 3.1 it is easy to show that
$\fp_i=\Z(S)\setminus V_i$ is a prime ideal. Now suppose that
$x\in V_i$ and $x^n=0$ and $x^{n-1}\neq 0$. Since $V_i\bigcup
\{0\}$ is an ideal of $S$ we have that $x^{n-1}\in V_i$. But
$x^n=x^{n-1} x=0$ and so $x^2=0$. We show that each part $V_i$
contains at most one nilpotent element. Let $x\neq y\in V_i$ are
two nilpotent elements. Then $xy\neq 0, y^2=x^2=0$. Therefore $xy$
is adjacent to $x$, which is a contradiction (note that $xy\in
V_i$). Now the assertion holds.
 Let $0\neq x\in S$ be a nilpotent element. By part (b), $x^2=0$.
There exists $1\le t\le r$ such that $x\in V_t$ and
$xSx=\{0\}\subseteq \fp_t$. Since $\fp_t$ is a prime ideal we have
that $x\in\fp_t$ and so $x=0$. This is a
contradiction.\hfill$\square$

\vspace{.2in}

In [{\bf AMY}, Theorem 3.5], it is shown that for an infinite ring
$R$, if $\Gamma(R)$ is a complete $r$-partite graph with $r\ge 3$
then $r$ is a power of a prime integer. The following example
shows that this is not true for semigroups. First we recall a
notion that we use in this example. Let $S_1, S_2 ,\cdots$ be
commutative semigroups with a zero element and $S_i\cap S_j=\{0\}$
whenever $i\neq j$, the 0-{\it orthogonal union} of $S_1,
S_2,\cdots$ is the semigroup $S=S_1\cup S_2\cup\cdots$ in which
every $S_i$ is a subsemigroup and $S_iS_j=0$ whenever $i\neq j$.

\vspace{.1in}

\noindent{\bf Example 3.7.} Let $S$ be the 0-orthogonal union of
$S_1,S_2,\cdots$. Let $|S_i|>2$ for all $i=1,2,\cdots ,r$. Then
$\Gamma(S)$ is a complete $r$-partite graph if and only if $\Z(S)$
is a 0-orthogonal union of semigroups without nonzero zero-divisors
(namely, the semigroups $S_i=V_i\cup\{0\}$).

\vspace{.2in}

\noindent{\bf Remark 3.8.} Note that in Example 3.7 the condition
$|S_i|>2$ is necessary. For example consider $S=\{0,a,b,c\}$ with
$ab=ac=a^2=0$, $bc=c^2=b^2=a$. In this case $\Gamma(S)$ is
complete bipartite and $S$ is not a 0-orthogonal union of non-zero
semigroups.

\vspace{.3in}

\noindent{\bf 4. Semigroups of finite chromatic number}

\vspace{.2in}

In this section, we begin to characterize the commutative
semigroups of finite chromatic number. Note that Beck in [{\bf B}]
and Anderson-Naseer in [{\bf AN}] let all elements of $R$ be
vertices of the graph $\Gamma(R))$ but we just consider the
nonzero zero-divisors. This is the reason why the chromatic number
(resp. clique number) of $\Gamma(R)$, in this paper, is one less
than the chromatic number (resp. clique number) of $\Gamma(R)$ in
[{\bf B}] and [{\bf AN}].

\vspace{.1in}

\noindent{\bf Theorem 4.1.} For a reduced semigroup $S$ the
following are equivalent:

\begin{itemize}

\item[(1)] $\chi(S)$ is finite.

\item[(2)] $\omega(S)$ is finite (i.e. $\Gamma(S)$ does not contain an infinite clique).

\item[(3)] The zero ideal in $S$ is a finite intersection of prime
ideals.

\end{itemize}

\vspace{.1in}

\noindent{\it Proof.} Since $\clique(S)\le\chi(S)$, the
implications (1)$\Rightarrow$(2) and is evident.

Now we prove (3)$\Rightarrow$(1). Let
$0=\fp_1\cap\fp_2\cap\cdots\cap\fp_k$ where for any $i$, $\fp_i$ is
a prime ideal. For any $0\neq x\in\Z(S)$, there exists minimum $j$,
such that $x\notin\fp_j$. Color $x$ with $j$. Now suppose that $x,y$
are colored to color $j$. If $xy=0$, then $xSy\subseteq\fp_j$. Since
$\fp_j$ is a prime ideal, then $x\in\fp_j$ or $y\in\fp_j$, which is
contradiction. So we have a $k$-coloring. Thus $\chi(S)\le k$.

It is now sufficient to show (2)$\Rightarrow$(3). By Proposition
2.7, $S$ satisfies the a.c.c. on annihilators. Let $T=\{\Ann
x_i|i\in I\}$ be the set of maximal members of the family $\{\Ann
a|a\neq 0\}$. By (4), Lemma (2.8) and Proposition 2.9(a), $T$ is a
finite set, and every element of $T$ is a prime ideal. Consider
$0\neq x\in S$. Then $\Ann x\subseteq\Ann x_i$ for some $i\in I$. If
$xx_i=0$, then $x_i\in\Ann x\subseteq\Ann x_i$, and so $x_i^2=0$.
Since $S$ is a reduced semigroup, then $x_i=0$, which is a
contradiction. Therefore $xx_i\neq 0$, and then $x\notin\Ann x_i$.
Thus $\cap_{i\in I}\Ann x_i=0$.\hfill$\square$

\vspace{.2in}

It is known that $\chi(G)\ge\omega(G)$ for general graph $G$ (see
[{\bf CO}, page 289]). Beck showed that if $R$ is a finite direct
product of reduced coloring and principal ideal rings then
$\chi(\Gamma(R))=\omega(\Gamma(R))$. In the following result the
equality $\chi(S)=\omega(S)$ is shown for some special case.

\vspace{.1in}

\noindent{\bf Corollary 4.2.} Suppose $S$ be a reduced semigroup.
Suppose $0=\cap_{i=1}^n\fp_i$ is a minimal prime decomposition of
$0$ (i.e. for any $i\neq j$, $\fp_i\neq\fp_j$ and for any $1\le t\le
n$, $0\neq\cap_{i\neq t}\fp_i$). Then $\chi(S)=\omega(S)=n$.

\vspace{.1in}

\noindent{\it Proof.} By the proof of Theorem 4.1, we have
$\chi(S)\le n$. Let $x_i\in\cap_{i\neq t}\fp_i\setminus\fp_t$. Then
$x_1,x_2,\cdots ,x_n$ is a clique and so $\omega(S)\ge n$. Now we
have $n\le\omega(S)\le\chi(S)\le n$, and hence
$\omega(S)=\chi(S)=n$.\hfill$\square$

\vspace{.2in}

\noindent{\bf Example 4.3.} Let $X$ be a $n$-set. We know that
$({\cal P}(X),\cap)$ is a reduced semigroup, where ${\cal P}(X)$
is the power set of $X$. For any $x\in X$, set $B_x=X-\{x\}$.
Clearly, for any $x\in X$, $({\cal P}(B_x),\cap)$ is a prime
ideal, and $\cap_{x\in X}{\cal P}(B_x)=\{\varnothing\}$. Thus
$\chi({\cal P}(X))=\omega({\cal P}(X))=n$.

\vspace{.2in}

Beck showed that for $n\le 3$, $\chi(\Gamma(R))=n$ if and only if
$\omega(\Gamma(R))=n$. Now we are ready to show that for $n\le 2$,
$\chi(\Gamma(S))=n$ if and only if $\omega(\Gamma(S)=n$.

\vspace{.1in}

\noindent{\bf Theorem 4.4.} Let $S$ be a commutative semigroup.
Then for $n\le 2$, $\chi(\Gamma(S))=n$ if and only if
$\omega(\Gamma(S)=n$.

\vspace{.1in}

\noindent{\it Proof.} The case $n=1$ is clear. If $\chi(S)=2$,
then $\Gamma(S)$ has at least two vertices and so $\omega(S)\ge
2$. On the other hand $\omega(S)\le\chi(S)=2$. Thus $\omega(S)=2$.

Conversely, let $\omega(S)=2$. If $\chi(S)>2$, then $\Gamma(S)$ is
not bipartite and so has a cycle of odd length. Let $C$ be the odd
cycle of minimal length. Since $\omega(S)=2$, the length of $C$ is
at least five (otherwise, the length of $C$ is $3$ and so
$\omega(S)=3$ that is a contradiction). Set
\begin{center}
$C$:$\,\,\,\, x_1$---$x_2$---$\cdots$---$x_n$---$x_1$,
\end{center}
where $n\ge 5$ is an odd integer. If $x_1x_3=0$, then $\Gamma(S)$
has a cycle of length $3$, which is a contradiction. Thus
$x_1x_3\neq 0$. Since all vertices in the cycle $C$ has degree $2$
and $x_1x_3$ has degree 3, we have $x_1x_3\neq x_i$ for any $1\le
i\le n$. Now consider the following cycle:
\begin{center}
$C'$:$\,\,\,\,
x_1x_3$---$x_4$---$x_5$---$\cdots$---$x_{n-1}$---$x_n$---$x_1x_3$.
\end{center}
It is easy to see that the length of $C'$ is $n-2$, which is a
contradiction. Thus $\Gamma(S)$ has no odd cycle. Therefore
$\Gamma(S)$ is bipartite and so $\chi(S)=2$.\hfill$\square$

\vspace{.2in}

Beck conjectured that $\chi(\Gamma(R))=\omega(\Gamma(R))$ in
general. In [{\bf AN}], Anderson and Naseer have given an example
of a finite local ring with $\chi(\Gamma(R))=5$ and
$\omega(\Gamma(R))=4$ thus giving a counterexample to Beck's
conjecture. For $n=1$ or $2$, $\chi(S)=n$ if and only if
$\omega(S)=n$. Now by giving an example we show that this result
is not true for $n=3$.

\vspace{.1in}

\noindent{\bf Example 4.5.} Let $S=\{0,a,b,c,d,e,f\}$ with
$fx=x^2=0$ for all $x\in S$. Also $ab=bc=cd=de=ae=0$, and
$ac=ad=bd=be=ce=f$. Then $\chi(S)=4$ and $\omega(S)=3$

\vspace{.3in}


\baselineskip=16pt

\begin{center}
\large {\bf References}
\end{center}
\vspace{.2in}

\begin{itemize}

\item[[AFL]] D. F. Anderson, A. Frazier, A. Lauve, P. S. Livingston,
{\it The Zero-Divisor Graph of a Commutative Ring II}, Lecture
Notes in Pure and Appl. Math., 220, Dekker, New York, 2001.

\item[[AL]] D. F. Anderson, P. S. Livingston, {\it The Zero-Divisor
Graph of a Commutative Ring}, J. Algebra {\bf 217} (1999), no. 2,
434-447.

\item[[ALS]] D. F. Anderson, R. Levy, J. Shapiro, {\it Zero-Divisor
Graphs, von Neumann Regular Rings, and Boolean Algebras}, J. Pure
Appl. Algebra {\bf 180} (2003), no. 3, 221-241.

\item[[AMY]] S. Akbari, H. R. Maimani, S. Yassemi, {\it When a
Zero-Divisor Graph is Planar or a Complete $r$-Partite Graph}, J.
Algebra {\bf 270} (2003), no. 1, 169-180.

\item[[AN]] D. D. Anderson, M. Naseer, {\it Beck's coloring of a commutative ring}, J. Algebra {\bf 159} (1991),
500--514.

\item[[B]] I. Beck, {\it Coloring of Commutative Rings}, J. Algebra
{\bf 116} (1988), no. 1, 208-226.

\item[[CO]] G. Chartrand, O. R. Oellermann, {\it Applied and
Algorithmic Graph Theory}, McGraw-Hill, Inc., New York, 1993.

\item[[DD]] F. DeMeyer, L. DeMeyer, {\it
Zero-Divisor Graphs of Semigroups}, J. Algebra {\bf 283} (2005),
190-198.

\item[[DMK]] F. R. DeMeyer, T. McKenzie, K. Schneider, {\it The
Zero-Divisor Graph of a Commutative Semigroup}, Semigroup Forum
{\bf 65} (2002), 206-214.

\item[[H]] J. M. Howie, {\it An introduction to semigroup theory}, L.M.S.
Monographs, No. 7. Academic Press [Harcourt Brace Jovanovich,
Publishers], London-New York, 1976.

\item[[LS]] R. Levy, J. Shapiro, {\it The Zero-Divisor Graph of von
Neumann Regular Rings}, Comm. Algebra {\bf 30} (2002), no. 2,
745-750.

\item[[ZW]] M. Zuo, T. Wu, {\it A New Graph Structure of Commutative
Semigroup}, Semigroup Forum {\bf 70} (2005), 71-80.

\end{itemize}

\end{document}